\documentclass{elsart}
\usepackage{amsmath, amssymb}

\renewcommand{\phi}{\varphi}

\newcommand{\BB}{\mathbb}
\newcommand{\cal}{\mathcal}
\newcommand{\g}{\mathfrak}
\newcommand{\separate}{\vskip11pt}
\newcommand{\tr}{\operatorname{Tr}}

\newtheorem{ex}[thm]{Example}

\journal{the Journal of Functional Analysis}

\begin{document}

\begin{frontmatter}

\title{\bf A fixed point localization formula for the
Fourier transform of regular semisimple coadjoint orbits}

\author{Matvei Libine}

\address{Department of Mathematics and Statistics,
University of Massachusetts,
710~North~Pleasant~Street, Amherst, MA 01003, USA}

\ead{matvei@math.umass.edu}

\tableofcontents

\begin{abstract}
Let $G_{\BB R}$ be a Lie group acting on an oriented manifold $M$,
and let $\omega$ be an equivariantly closed form on $M$.
If both $G_{\BB R}$ and $M$ are compact, then the integral $\int_M \omega$
is given by the fixed point integral localization formula
(Theorem 7.11 in \cite{BGV}).
Unfortunately, this formula fails when the acting Lie group $G_{\BB R}$
is not compact: there simply may not be enough fixed points present.
A proposed remedy is to modify the action of $G_{\BB R}$ in such a way
that all fixed points are accounted for.

Let $G_{\BB R}$ be a real semisimple Lie group, possibly {\em noncompact}.
One of the most important examples of equivariantly closed forms is the
symplectic volume form $d\beta$ of a coadjoint orbit $\Omega$.
Even if $\Omega$ is not compact, the integral $\int_{\Omega} d\beta$
exists as a distribution on the Lie algebra $\g g_{\BB R}$.
This distribution is called
{\em the Fourier transform of the coadjoint orbit}.

In this article we will apply the localization results described in
\cite{L1} and \cite{L2} to get a geometric derivation of
Harish-Chandra's formula (\ref{charformula})
for the Fourier transforms of regular semisimple coadjoint orbits.
Then we will make an explicit computation for
%semisimple coadjoint orbits of $SL(2,\BB R)$ and
the coadjoint orbits of elements of $\g g_{\BB R}^*$
which are dual to regular semisimple elements lying in
a maximally split Cartan subalgebra of $\g g_{\BB R}$.
\end{abstract}

\begin{keyword}
equivariant forms, fixed point integral localization formula,
semisimple Lie groups, Fourier transform of a coadjoint orbit,
invariant eigendistributions, integral character formula,
characteristic cycles of sheaves
\end{keyword}

\end{frontmatter}

%\separate

\begin{section}
{Introduction}
\end{section}

Equivariant forms were introduced in 1950 by Henri Cartan.
There are many good texts on this subject including \cite{BGV}
and \cite{GS}.
Let $G_{\BB R}$ be a real Lie group acting on a compact oriented manifold $M$,
let $\g g_{\BB R}$ be the Lie algebra of $G_{\BB R}$, and let
$\omega : \g g_{\BB R} \to {\cal A}(M)$ be an equivariantly closed form on $M$.
(Here ${\cal A}(M)$ denotes the algebra of smooth differential forms on $M$.)
For $X \in \g g_{\BB R}$, we denote by $M^X$ the set of zeroes of the vector
field on $M$ induced by the infinitesimal action of $X$.
We assume that $M^X$ is discrete.
Then the integral localization formula (Theorem 7.11 in \cite{BGV}) says
that when the Lie group $G_{\BB R}$ is {\em compact} the integral of
$\omega(X)$ can be expressed as a sum over the
set of zeroes $M^X$ of certain {\em local} quantities of $M$ and $\omega$:
$$
\int_M \omega(X) = \sum_{p \in M^X}
\text{local invariant of $M$ and $\omega$ at $p$}.
$$
This is the essence of the integral localization formula for equivariantly
closed differential forms.

Unfortunately, this formula fails when the acting Lie group $G_{\BB R}$
is not compact: there simply may not be enough fixed points present.
For example, let $G_{\BB R} = SL(2, \BB R)$. Then there exists an
essentially unique $G_{\BB R}$-invariant symmetric bilinear form $B$ on
$\g g_{\BB R} = \g {sl}(2,\BB R)$, say $B(X,Y)$ is the Killing form
$\tr(ad(X)ad(Y))$.
This bilinear form $B$ induces an $SL(2,\BB R)$-equivariant isomorphism
$I: \g {sl}(2,\BB R) \, \tilde \to \, \g {sl}(2,\BB R)^*$.
Let $\Omega \subset \g {sl}(2,\BB R)^*$ denote the coadjoint orbit of
$I \begin{pmatrix} 0 & 1 \\ -1 & 0 \end{pmatrix}$.
Like all semisimple coadjoint orbits, $\Omega$ possesses a canonical symplectic
form $\omega$ which is the top degree part of a certain equivariantly
closed form. Although $\Omega$ is not compact, the symplectic volume
$\int_{\Omega} \omega$ still exists as a distribution on $\g {sl}(2,\BB R)$.
Let $\g {sl} (2,\BB R)'_{split}$ denote the set of
$X \in \g {sl}(2,\BB R)$ with distinct real eigenvalues.
The set $\g {sl} (2,\BB R)'_{split}$ is an open subset of $\g {sl} (2,\BB R)$.
Now, if we take any element $X \in \g {sl} (2,\BB R)'_{split}$,
then one can see that the vector field on $\Omega$ generated by $X$
has no zeroes. Thus if there were a fixed point integral localization
formula like in the case of compact group $G_{\BB R}$,
this formula would suggest that the distribution
determined by $\int_{\Omega} \omega$ vanishes on the open set
$\g {sl} (2,\BB R)'_{split}$.
But it is known that the restriction $\int_{\Omega} \omega$ to
$\g {sl} (2,\BB R)'_{split}$ is {\em not} zero.
%We will see in Section \ref{SU(1,1)} that it is not so.

On the other hand, the localization argument for character formulas
described in \cite{L1} and \cite{L2} strongly suggests that there
should be some kind of a localization formula when the group $G_{\BB R}$ 
is semisimple. One can regard the failure of the localization formula
in the last example as an instance when some of the fixed points which
should be counted are ``hiding somewhere at infinity''.
A proposed remedy is to modify the action of $G_{\BB R}$ in such a way
that all fixed points are accounted for.

We will show how it can be done in the case of the symplectic volume of
a regular semisimple coadjoint orbit.
Let $\lambda \in \g g_{\BB R}^*$ be a regular semisimple element and
let $\Omega_{\lambda, \BB R}$ be its coadjoint orbit.
Let ${\cal X}$ denote the {\em flag variety}, i.e the variety
of Borel subalgebras
$\g b \subset \g g = \g g_{\BB R} \otimes_{\BB R} \BB C$.
The space ${\cal X}$ is a smooth complex projective variety.
The twisted moment map $\mu_{\lambda}: T^*{\cal X} \to \g g^*$
defined in \cite{SchV1} is an embedding whose image contains
$\Omega_{\lambda, \BB R}$. We will replace $\Omega_{\lambda, \BB R}$ with
$\mu_{\lambda}^{-1}(\Omega_{\lambda, \BB R})$ and then it will be possible to
apply the results of \cite{SchV2} and
the localization argument described in \cite{L1} and \cite{L2}
to get a fixed point integral formula (\ref{charformula})
for the Fourier transform of $\Omega_{\lambda, \BB R}$.
This is not really a new formula, but rather a different way of looking at
Harish-Chandra's formula for an invariant eigendistribution.
Unless the group $G_{\BB R}$ is compact, the map
$\mu_{\lambda}: T^*{\cal X} \to \g g^*$ is not $G_{\BB R}$-equivariant.
Thus the cycle $\mu_{\lambda}^{-1}(\Omega_{\lambda, \BB R})$ does not
need to be $G_{\BB R}$-invariant, but we will show that it is homologous
to a $G_{\BB R}$-invariant cycle which we will call $C_{\lambda}$.
Replacing $\Omega_{\lambda, \BB R}$ with $C_{\lambda}$
amounts to changing the action of $G_{\BB R}$ on $\Omega_{\lambda, \BB R}$.

As an example, we will compute the Fourier transforms of the coadjoint
orbits of elements of $\g g_{\BB R}^*$ which are dual to regular semisimple
elements lying in a maximally split Cartan subalgebra of $\g g_{\BB R}$.

In this setting the main difficulty lies in identifying the cycle
$C_{\lambda}$.
It is most likely that this cycle $C_{\lambda}$ can be identified as
the characteristic cycle of some sheaf ${\cal F}_{\lambda}$ which has
a simple description.
Then one will get an interesting perspective of looking at Harish-Chandra's
formula for the Fourier transforms of regular semisimple coadjoint orbits.

%In this article we will compute the Fourier transforms of all regular
%semisimple coadjoint orbits of $SL(2,\BB R)$ and
%the coadjoint orbits of elements of
%$\g g_{\BB R}^*$ which are dual to regular semisimple elements lying in
%a maximally split Cartan subalgebra of $\g g_{\BB R}$.

This result strongly resembles the following result due to W.~Rossmann
which is described in \cite{BGV}, \cite{BV}, \cite{Par}, \cite{R2}, \cite{V}.
Let $G_{\BB R}$ be a connected real semisimple Lie group, and let
$K_{\BB R} \subset G_{\BB R}$ be a maximal compact subgroup of
$G_{\BB R}$ with Lie algebra $\g k_{\BB R}$.
Let $\g t_{\BB R} \subset \g k_{\BB R}$ be a Cartan subalgebra of
$\g k_{\BB R}$ and suppose that $\g t_{\BB R}$ is also a Cartan subalgebra of
$\g g_{\BB R}$. Then W.~Rossmann gives
an explicit formula for the restriction of the Fourier transform of
a semisimple coadjoint orbit to $\g t_{\BB R}$.

The difference between the result described in this article and that
of W.~Rossmann is that it seldom happens that a Cartan
subalgebra $\g t_{\BB R} \subset \g k_{\BB R}$ is also a Cartan subalgebra of
$\g g_{\BB R}$. Also, W.~Rossmann does not give an explicit formula 
for the Fourier transform as a distribution on all of $\g g_{\BB R}$, but
rather the restriction of the Fourier transform to $\g t_{\BB R}$.

I would like to thank W.~Rossmann for his helpful suggestions and
especially for pointing out that part {\it (1)} of Theorem \ref{main}
was stated and used in his proof of Theorem 1.4(b) in \cite{R3}.
(The proof of this statement was omitted because it was not the point
of that paper.)

This article strongly suggests that some results which previously
were known in the compact group setting only can be generalized to
noncompact groups. It would be very interesting to have a uniform
equivariant form theory which holds for both compact and noncompact
groups and which generalizes presently known compact case.

%\separate

\begin{section}
{The Fourier Transform of a Coadjoint Orbit}
\end{section}  

We fix a {\em connected, complex algebraic, semisimple} Lie group $G$
which is defined over $\BB R$.
Let $G_{\BB R}$ be the identity component of the group of real points
$G(\BB R)$.
We denote by $\g g$ and $\g g_{\BB R}$ the Lie algebras of
$G$ and $G_{\BB R}$ respectively.
The group $G$ (respectively $G_{\BB R}$) acts on $\g g$
(respectively $\g g_{\BB R}$) by the {\em adjoint} action, and the dual
action of  $G$ (respectively $G_{\BB R}$) on $\g g^*$
(respectively $\g g_{\BB R}^*$) is called the {\em coadjoint} action.

Let $B(X,Y)$ be an (essentially unique) $G$-invariant symmetric
bilinear form on $\g g$, say $B(X,Y)$ is the Killing form $\tr(ad(X)ad(Y))$.
This bilinear form $B$ induces a $G$-equivariant isomorphism
$I: \g g \, \tilde \to \, \g g^*$ which restricts to a
$G_{\BB R}$-equivariant isomorphism
$\g g_{\BB R} \, \tilde \to \, \g g_{\BB R}^*$.

We denote by $\g g_{\BB R}'$ the set of {\em regular semisimple} elements
in $\g g_{\BB R}$.
These are elements $X \in \g g_{\BB R}$ such that the action of
$ad(X)$ on $\g g_{\BB R}$ can be diagonalized over $\BB C$ and $ad(X)$
has maximal possible rank.
We also call the $I$-image of $\g g_{\BB R}'$
the {\em regular semisimple} elements of $\g g_{\BB R}^*$ and
denote them by $(\g g_{\BB R}^*)' = I (\g g_{\BB R}')$.

The set $\g g_{\BB R}'$ is an open subset of $\g g_{\BB R}$ and its
complement has measure zero.
Each regular semisimple $X \in \g g'_{\BB R}$ lies in a unique Cartan
subalgebra $\g t_{\BB R}(X) \subset \g g_{\BB R}$.
It is important to remember that when $G_{\BB R}$ is not compact some Cartan
subalgebras of $\g g_{\BB R}$ may not be conjugate by an element from
$G_{\BB R}$.

\separate

We pick a regular semisimple element $\lambda' \in (\g g_{\BB R}^*)'$
and consider its coadjoint orbit
$G_{\BB R} \cdot \lambda' \subset \g g_{\BB R}^*$.
Because we are going to apply results of \cite{L1} and \cite{L2}
it will be more convenient to replace $\lambda'$ with
$\lambda = i \lambda' \in i(\g g_{\BB R}^*)'$ and consider
$$
\Omega_{\lambda,\BB R} = G_{\BB R} \cdot \lambda
\subset i \g g_{\BB R}^* \subset \g g^*.
$$
The $G_{\BB R}$-orbit $\Omega_{\lambda,\BB R}$ is a real submanifold of
the complex coadjoint orbit
$\Omega_{\lambda} = G \cdot \lambda \subset \g g^*$.

If $X \in \g g$ is an element of the Lie algebra,
we denote by $X_{\Omega_{\lambda}}$ the vector field on $\Omega_{\lambda}$
generated by $X$: if $\nu \in \Omega_{\lambda}$ and
$f \in {\cal C}^{\infty}(\Omega_{\lambda})$, then
$$
X_{\Omega_{\lambda}} (\nu)f=
\frac d{d\epsilon} f \bigl( \exp(\epsilon X) \cdot \nu \bigr)
\bigl|_{\epsilon=0}.
$$
We call a point $\nu \in \Omega_{\lambda}$ a {\em fixed point of $X$}
if the vector field $X_{\Omega_{\lambda}}$ on $\Omega_{\lambda}$
vanishes at $\nu$, i.e. $X_{\Omega_{\lambda}} (\nu)= 0$.

The complex coadjoint orbit $\Omega_{\lambda}$ has a canonical structure of
a complex symplectic manifold.
We describe its symplectic form $\sigma_{\lambda}$.
The vector fields $X_{\Omega_{\lambda}}$, $X \in \g g$, are sections of
the tangent bundle $T \Omega_{\lambda}$
and span the tangent space at each point.
Let $\nu \in \Omega_{\lambda}$, we define
$\sigma_{\lambda}$ on $T_{\nu} \Omega_{\lambda}$ by
$$
\sigma_{\lambda}(X_{\Omega_{\lambda}}, Y_{\Omega_{\lambda}})_{\nu} =_{def}
\langle \nu, [X,Y] \rangle, \qquad X,Y \in \g g.
$$
Then $\sigma_{\lambda}$ is a
well-defined $G$-invariant symplectic form on $\Omega_{\lambda}$,
furthermore the action of $G$
on $\Omega_{\lambda}$ is Hamiltonian and the symplectic moment of $X$
is the function $\nu \mapsto \langle \nu,X \rangle$.

We define the {\em Liouville form} on $\Omega_{\lambda}$ by
(notice the $\frac 1{(2\pi i)^n}$ factor)
\begin{equation}  \label{beta}
d\beta = \frac 1{(2\pi i)^n n!}  \sigma_{\lambda}^n,
\qquad n = \dim_{\BB C} (\Omega_{\lambda})/2
= \dim_{\BB R} (\Omega_{\lambda, \BB R})/2.
\end{equation}
The form $d\beta$ restricts to a volume form on $\Omega_{\lambda, \BB R}$
which determines a canonical orientation on $\Omega_{\lambda, \BB R}$.

\separate

At this point we need to specify that if $M$ is a complex manifold and
$M^{\BB R}$ is the underlying real analytic manifold,
then there are at least two different but equally natural identifications
between the holomorphic cotangent bundle $T^*M$ and the ${\cal C}^{\infty}$
cotangent bundle $T^*M^{\BB R}$; we use the convention of
\cite{KaScha}, Chapter XI. The same convention is used in
\cite{L1}, \cite{L2}, \cite{R1} and \cite{SchV2}. Under this convention, if
$\sigma_M$ is the canonical complex symplectic form on $T^*M$ and
$\sigma_{M^{\BB R}}$ is the canonical real symplectic form on $T^*M^{\BB R}$,
then $\sigma_{M^{\BB R}}$ gets identified with
$2 \operatorname{Re} \sigma_M$.

\separate

From now on we regard $\Omega_{\lambda, \BB R}$ as a Borel-Moore cycle in
$\Omega_{\lambda}$. We will use notation $|\Omega_{\lambda, \BB R}|$
to denote the support of $\Omega_{\lambda, \BB R}$.

\separate

If $\phi$ is a smooth compactly supported differential
form on $\g g_{\BB R}$ of top degree, then we define its Fourier transform
as in \cite{L1}, \cite{L2} and \cite{SchV2}:
$$
\hat \phi(\zeta) = \int_{\g g_{\BB R}} e^{\langle g, \zeta \rangle} \phi
\qquad (g \in \g g_{\BB R},\:\zeta \in \g g^*),
$$
without the customary factor of $i=\sqrt{-1}$ in the exponent.
Because $\hat \phi$ decays fast in the imaginary directions, the integral
$$
\widehat{\Omega_{\lambda, \BB R}}(\phi) =_{def}
\int_{\Omega_{\lambda, \BB R}} \hat \phi d\beta =
\frac 1{(2\pi i)^n n!}
\int_{\Omega_{\lambda, \BB R}} \hat \phi \sigma_{\lambda}^n
$$
is finite.
Thus we obtain a distribution on $\g g_{\BB R}$ called the
{\em Fourier transform} of the coadjoint orbit $\Omega_{\lambda, \BB R}$.
Since the symplectic form $\sigma_{\lambda}$ is $G_{\BB R}$-invariant,
so is the distribution $\widehat{\Omega_{\lambda, \BB R}}$.

\separate

Let ${\cal Z}({\cal U} (\g g_{\BB R}))$ denote the center of the
universal enveloping algebra of $\g g_{\BB R}$.
It is canonically isomorphic to the algebra of conjugate-invariant
constant coefficient differential operators on $\g g_{\BB R}$.
It is not hard to see that $\widehat{\Omega_{\lambda, \BB R}}$
is an eigendistribution with respect to
${\cal Z}({\cal U} (\g g_{\BB R}))$, i.e. each element of
${\cal Z}({\cal U} (\g g_{\BB R}))$ acts on 
$\widehat{\Omega_{\lambda, \BB R}}$ by multiplication by some scalar.
Then by Harish-Chandra's regularity theorem
(\cite{HC} or Theorem 3.3 in \cite{A}),
\begin{equation}  \label{F}
\widehat{\Omega_{\lambda, \BB R}} (\phi) =
\int_{\g g_{\BB R}} F_{\lambda} \phi,
\end{equation}
where $F_{\lambda}$ is an $Ad(G_{\BB R})$-invariant, locally $L^1$ function
on $\g g_{\BB R}$ which is represented by a real analytic function on
$\g g_{\BB R}'$.

%\separate

\begin{section}
{Harish-Chandra's Formula for $\widehat{\Omega_{\lambda, \BB R}}$}
\end{section}  

In this section we derive Harish-Chandra's formula (\ref{charformula}) for
the Fourier transform of a regular semisimple coadjoint orbit
$\widehat{\Omega_{\lambda, \BB R}}$ as a fixed point
integral localization formula.

\separate

Let ${\cal X}$ denote the {\em flag variety}, i.e the variety
of Borel subalgebras
$\g b \subset \g g = \g g_{\BB R} \otimes_{\BB R} \BB C$.
The space ${\cal X}$ is a smooth complex projective variety.

As before, if $X \in \g g$ is an element of the Lie algebra,
we denote by $X_{{\cal X}}$ the vector field on ${\cal X}$
generated by $X$: if $x \in {\cal X}$ and
$f \in {\cal C}^{\infty}({\cal X})$, then
$$
X_{{\cal X}} (x)f=
\frac d{d\epsilon} f \bigl( \exp(\epsilon X) \cdot x \bigr)
\bigl|_{\epsilon=0}.
$$

The ordinary moment map $\mu: T^*{\cal X} \to \g g^*$ is defined by
$$
\mu(\zeta): X \mapsto \langle \zeta, X_{\cal X} \rangle,
\qquad \zeta \in T^*{\cal X},\: X \in \g g^*.
$$
The moment map $\mu$ takes values in the nilpotent cone
${\cal N}^* \subset \g g^*$.

\separate

Next we define Rossmann's {\em twisted moment map}. It is best described
in \cite{SchV1}, Section 8.
Let $\g h$ be the {\em universal Cartan algebra} of $\g g$. That is
$\g h$ is a Cartan algebra equipped with a choice of positive roots.
Let $\tilde \lambda \in \g h^*$, and fix a compact real form
$U_{\BB R} \subset G$. The twisted moment map
$\mu_{\tilde \lambda}: T^*{\cal X} \to \g g^*$ can be written as
$$
\mu_{\tilde \lambda} = \mu + \tilde \lambda_x,
$$
where $\mu: T^*{\cal X} \to \g g^*$ is the ordinary moment map and
the map $\tilde \lambda_x: {\cal X} \to \g g^*$ is defined as follows.
For each $x \in {\cal X}$, let $B(x) \subset G$ and
$T_{\BB R}(x) \subset U_{\BB R}$ be the stabilizer groups of $x$:
\begin{equation}  \label{B(x)}
B(x) = \{ g \in G;\: g \cdot x = x \}, \qquad
T_{\BB R}(x) = \{ u \in U_{\BB R};\: u \cdot x = x \}.
\end{equation}
Then $B(x)$ is the Borel subgroup of $G$ corresponding to $x \in {\cal X}$,
and $T_{\BB R}(x)$ is the maximal torus of $U_{\BB R}$ fixing $x$.
Let $\g b_x$ denote the Lie algebra of $B(x)$. The quotient
$\g b_x/[\g b_x,\g b_x]$ is canonically isomorphic to the universal
Cartan algebra $\g h$ so that the nilpotent algebra
$\g n_x = [\g b_x,\g b_x]$ is the vector space
direct sum of the negative root spaces.
On the other hand, the complexified Lie algebra $\g t_x$
of $T_{\BB R}(x)$ can be realized as a subalgebra of $\g b_x$ which
maps isomorphically onto $\g b_x/[\g b_x,\g b_x]$.
Thus we obtain an identification of $\g t_x$ with the universal Cartan
algebra $\g h$. This identification makes $\tilde \lambda \in \g h^*$
correspond to a $\tilde \lambda_x \in (\g t_x)^*$, which we extend to a
linear functional on $\g g$ using the canonical splitting of
vector spaces
$$
\g g = \g t_x \oplus [\g t_x, \g g]
=\g t_x \oplus \bigl( \bigoplus \text{root spaces} \bigr).
$$
The map $x \mapsto \tilde \lambda_x \in \g g^*$ is real algebraic and
$U_{\BB R}$-equivariant (but not $G$-equivariant unless $\tilde \lambda =0$).

If $\tilde \lambda=0$, then $\tilde \lambda_x \equiv 0$ and
the twisted moment map $\mu_{\tilde \lambda} = \mu + \tilde \lambda_x$
reduces to the ordinary moment map $\mu$.
On the other hand, if $\tilde \lambda$ is regular, then
$\mu_{\tilde \lambda}$ takes values in the coadjoint orbit
$\Omega_{\tilde \lambda}$ and
$$
\mu_{\tilde \lambda}: T^*{\cal X} \, \tilde \longrightarrow \,
\Omega_{\tilde \lambda}
$$
is a real algebraic diffeomorphism which is equivariant with respect to
$U_{\BB R}$, but never $G$-equivariant (see Section 8 of \cite{SchV1}).

\separate

Recall that $\lambda \in i(\g g_{\BB R}^*)'$ is a regular semisimple element
and we are interested in the Fourier transform of
$\Omega_{\lambda,\BB R} = G_{\BB R} \cdot \lambda$
which is a cycle in $\Omega_{\lambda} = G \cdot \lambda$
such that its support $|\Omega_{\lambda,\BB R}|$ is a smooth submanifold.
Let $\g t$ be the unique Cartan subalgebra of $\g g$ which contains
$I^{-1} (i\lambda) \in \g g_{\BB R}' \subset \g g$, and let us select
some positive system of roots $\Phi^+$ on $\g t$.
Then we get an identification of $\g t$ with the universal Cartan
algebra $\g h$, and $\lambda$ can be regarded as a regular element
of $\g h^*$ and thus we obtain a twisted moment map
$\mu_{\lambda}: T^*{\cal X} \, \tilde \to \, \Omega_{\lambda}$.
Because $\mu_{\lambda}$ is a diffeomorphism, we can rewrite the Fourier
transform of $\Omega_{\lambda, \BB R}$ as
\begin{equation}  \label{FT1}
\widehat{\Omega_{\lambda, \BB R}}(\phi) =
\frac 1{(2\pi i)^n n!}
\int_{\Omega_{\lambda, \BB R}} \hat \phi \sigma_{\lambda}^n =
\frac 1{(2\pi i)^n n!}
\int_{\mu_{\lambda}^{-1}(\Omega_{\lambda, \BB R})}
\mu_{\lambda}^* \hat \phi \cdot (\mu_{\lambda}^* \sigma_{\lambda})^n.
\end{equation}

The form $\mu_{\lambda}^* \sigma_{\lambda}$ is known -- it is
precisely the form that appears in the integral character formula
in \cite{L1}, \cite{L2}, \cite{R1} and \cite{SchV2}.
In order to make a precise statement, let $\tau_{\lambda}$ be a
$U_{\BB R}$-invariant real algebraic two form on ${\cal X}$ defined
by the formula
$$
\tau_{\lambda} (u_x, v_x) = \lambda_x([u,v]),
$$
where $u_x, v_x \in T_x {\cal X}$ are the tangent vectors at $x$ induced by
$u, v \in \g u_{\BB R}$ via the infinitesimal action of the Lie algebra
$\g u_{\BB R}$ of $U_{\BB R}$.
Let $\sigma$ be the canonical complex algebraic $G$-invariant
symplectic form on $T^*{\cal X}$,
and let $\pi: T^*{\cal X} \to {\cal X}$ be the natural projection.
Then Proposition 3.3 of \cite{SchV2} (which was derived from
Lemma 7.2 in \cite{R1}) says that
\begin{equation}  \label{pullback}
\mu_{\lambda}^* \sigma_{\lambda} = - \sigma + \pi^* \tau_{\lambda}.
\end{equation}
This form appears in the integral character formula
(Theorem 2.1 in \cite{R1}, Theorem 3.8 in \cite{SchV2},
formula (2) in \cite{L1} and formula (1) in \cite{L2}).

It is a much more difficult task to identify the cycle
$\mu_{\lambda}^{-1}(\Omega_{\lambda, \BB R})$ in $T^*{\cal X}$.
Because $\mu_{\lambda}$ is not $G_{\BB R}$-equivariant, this cycle
even may not be $G_{\BB R}$-invariant. But, as we will see soon,
$\mu_{\lambda}^{-1}(\Omega_{\lambda, \BB R})$ is
homologous to a conic Lagrangian $G_{\BB R}$-invariant cycle.

\separate

\begin{lem}  \label{bounded}
The composition map
$\mu \circ \mu_{\lambda}^{-1} : \Omega_{\lambda} \to \g g^* \simeq
\g g_{\BB R}^* \oplus i \g g_{\BB R}^*$ sends
$|\Omega_{\lambda, \BB R}|$ into a set with bounded real part, i.e.
the set
$$
\operatorname{Re}
\bigl( \mu \circ \mu_{\lambda}^{-1} (|\Omega_{\lambda, \BB R}|) \bigr)
$$
is a bounded subset of $\g g_{\BB R}^*$.
\end{lem}

\pf
Notice that $\mu = \mu_{\lambda} - \lambda_x$. Hence
$\mu \circ \mu_{\lambda}^{-1} =
j_{\Omega_{\lambda} \hookrightarrow \g g^*} -
\lambda_x \circ \pi \circ \mu_{\lambda}^{-1}$,
where $j_{\Omega_{\lambda} \hookrightarrow \g g^*}$
denotes the inclusion map.
Since $|\Omega_{\lambda, \BB R}| \subset i \g g_{\BB R}^*$,
$$
\operatorname{Re}
\bigl( \mu \circ \mu_{\lambda}^{-1} (|\Omega_{\lambda, \BB R}|) \bigr) =
\operatorname{Re}
 \bigl( -\lambda_x \circ \pi \circ \mu_{\lambda}^{-1}
(|\Omega_{\lambda, \BB R}|) \bigr)
\subset \operatorname{Re} (-\lambda_x ({\cal X})).
$$
Because the flag variety ${\cal X}$ is compact and the map
$\lambda_x$ is continuous the result follows.
\qed

If $s$ is a positive real number, then the fiber scaling map
$$
Scaling_s: T^*{\cal X} \to T^*{\cal X}, \qquad \zeta \mapsto s\zeta,
$$
sends
the cycle $\mu_{\lambda}^{-1}(\Omega_{\lambda, \BB R})$ into
a new cycle denoted by $s \cdot \mu_{\lambda}^{-1}(\Omega_{\lambda, \BB R})$.
We use the notion of families of cycles and their limits described in
Section 3 of \cite{SchV1}.
The cycles $s \cdot \mu_{\lambda}^{-1}(\Omega_{\lambda, \BB R})$,
$s \in (0,1)$, piece together to form a family of cycles $C_{(0,1)}$
parameterized by an open interval $(0,1)$.
Thus we can take its limit as $s \to 0^+$: let
\begin{equation}  \label{C}
C_{\lambda} =
\lim_{s \to 0^+} s \cdot \mu_{\lambda}^{-1}(\Omega_{\lambda, \BB R}).
\end{equation}
It is clear from the construction that the cycle $C_{\lambda}$ is {\em conic},
i.e. invariant under the fiber scaling maps
$Scaling_s: T^*{\cal X} \to T^*{\cal X}$, $s>0$.

In general, if $\tilde C_{(0,a)}$ is a family of cycles
in some space $Z$ parameterized by an open interval $(0,a)$, then
$\lim_{s \to 0^+} \tilde C_s = \tilde C_0$
means that if we regard $\tilde C_{(0,a)}$ as a chain in $[0,a) \times Z$,
then $\partial \tilde C_{(0,a)} = - \tilde C_0$.
Thus it is natural to define $\lim_{s \to a^-} \tilde C_s$
by setting
$$
\lim_{s \to a^-} \tilde C_s =
\partial \tilde C_{(0,a)},
$$
where $\tilde C_{(0,a)}$ is regarded as a chain in $(0,a] \times Z$.

\separate

When $s=1$, the fiber scaling map
$Scaling_1: T^*{\cal X} \to T^*{\cal X}$ is the identity map
on $T^*{\cal X}$, and it is clear that
$\lim_{s \to 1^-} s \cdot \mu_{\lambda}^{-1}(\Omega_{\lambda, \BB R}) =
\mu_{\lambda}^{-1}(\Omega_{\lambda, \BB R})$.
Next we make the following observation.
Recall that the family of cycles $C_{(0,1)}$ is a cycle in
$(0,1) \times T^*{\cal X}$.
We can regard it as a chain in $[0,1] \times T^*{\cal X}$.
Let $p: [0,1] \times T^*{\cal X} \twoheadrightarrow T^*{\cal X}$
be the projection.
Then $p_*(C_{(0,1)})$ is a Borel-Moore
$(2n+1)$-chain in $T^*{\cal X}$ such that
$$
\partial \bigl( p_*(C_{(0,1)}) \bigr)
= \lim_{s\to 1^-} C_s - \lim_{\epsilon \to 0^+} C_s
= \mu_{\lambda}^{-1}(\Omega_{\lambda, \BB R}) - C_{\lambda}.
$$
In particular, the cycles $C_{\lambda}$
and $\mu_{\lambda}^{-1}(\Omega_{\lambda, \BB R})$ are homologous.
Moreover, the support
of $p_*(C_{(0,1)})$ lies in the closure
$$
\overline{
\bigcup_{s \in [0,1]} s|\mu_{\lambda}^{-1}(\Omega_{\lambda, \BB R})|}.
$$

We see from Lemma \ref{bounded} that the image under
$\operatorname{Re} (\mu): T^*{\cal X} \to \g g_{\BB R}^*$ of the support
$|p_*(C_{(0,1)})|$ is bounded and hence Lemma 3.19 from \cite{SchV2}
together with equation (\ref{pullback}) implies that
\begin{equation}  \label{FT2}
\int_{\mu_{\lambda}^{-1}(\Omega_{\lambda, \BB R})}
\mu_{\lambda}^* \hat \phi \cdot (\mu_{\lambda}^* \sigma_{\lambda})^n
= \int_{C_{\lambda}}
\mu_{\lambda}^* \hat \phi \cdot (\mu_{\lambda}^* \sigma_{\lambda})^n.
\end{equation}

\separate

The group $G_{\BB R}$ acts on ${\cal X}$ real algebraically with
finitely many orbits. The resulting stratification of ${\cal X}$ by
$G_{\BB R}$-orbits satisfies the Whitney conditions. We denote by
$T^*_{G_{\BB R}} {\cal X}$ the union of the ${\cal C}^{\infty}$
conormal bundles of the $G_{\BB R}$-action.
The set $T^*_{G_{\BB R}} {\cal X}$ is a closed $G_{\BB R}$-invariant
real semialgebraic subset of $T^*{\cal X}$ of dimension $2n$.
We can also express this set as
$$
T^*_{G_{\BB R}} {\cal X} = \mu^{-1} (i \g g_{\BB R}^*).
$$
Clearly, the moment map $\mu$ takes only purely imaginary values on
the support $|C_{\lambda}|$. Thus $|C_{\lambda}|$ is a subset of
$T^*_{G_{\BB R}} {\cal X}$.
This makes $C_{\lambda}$ a top-dimensional cycle in $T^*_{G_{\BB R}} {\cal X}$.
Because the group $G_{\BB R}$ was assumed to be connected, it follows that
the cycle $C_{\lambda}$ is $G_{\BB R}$-invariant.

\separate

We already know that $C_{\lambda}$ is a conic cycle, and since its support
$|C_{\lambda}|$ lies in a finite union of conormal bundles,
$C_{\lambda}$ is a Lagrangian
cycle with respect to the canonical symplectic structure on
$T^*{\cal X}$, viewed as a ${\cal C}^{\infty}$ (not holomorphic!)
cotangent bundle on ${\cal X}$. Hence by Theorem 9.7.1 of \cite{KaScha},
there exists a bounded complex of sheaves ${\cal F}_{\lambda}$
on ${\cal X}$ with $\BB R$-constructible cohomology such that its
characteristic cycle $Ch({\cal F}_{\lambda}) = C_{\lambda}$.

Characteristic cycles were introduced by M.~Kashiwara and their
definition can be found in \cite{KaScha}. On the other hand,
W.~Schmid and K.~Vilonen give a geometric way to understand
characteristic cycles in \cite{SchV1}.
A comprehensive treatment of characteristic cycles can be found in
\cite{Schu}.

Combining this observation with equations (\ref{FT1}), (\ref{pullback})
and (\ref{FT2}) we can rewrite the Fourier transform of
$\Omega_{\lambda, \BB R}$ as follows
\begin{equation}  \label{FT}
\widehat{\Omega_{\lambda, \BB R}}(\phi) =
\frac 1{(2\pi i)^n n!} \int_{Ch({\cal F}_{\lambda})}
\mu_{\lambda}^* \hat \phi \cdot (-\sigma + \pi^* \tau_{\lambda})^n.
\end{equation}
The right hand side of this formula is precisely the integral character
formula from \cite{L1}, \cite{L2} and \cite{SchV2}.
If we can identify the complex of sheaves ${\cal F}_{\lambda}$,
then we can apply the localization argument from \cite{L1}
to compute this integral and identify the function $F_{\lambda}$ in
the equation (\ref{F}).
(We do not require the complex of sheaves ${\cal F}_{\lambda}$ to be
$G_{\BB R}$-equivariant or to be an element of the
``$G_{\BB R}$-equivariant derived category on ${\cal X}$ with twist
$(\lambda -\rho)$,'' $\operatorname{D}_{G_{\BB R}}({\cal X})_{\lambda}$,
as in \cite{SchV2} because the
localization argument in \cite{L1} only uses that
$Ch({\cal F}_{\lambda})$ is $G_{\BB R}$-invariant.)

If $X \in \g g_{\BB R}'$ is a regular semisimple element,
let ${\cal X}^X$ denote the set of zeroes of the vector field
$X_{\cal X}$ on ${\cal X}$, and let $\g t(X) \subset \g g$
be the unique Cartan subalgebra containing $X$. Then
\begin{equation}  \label{charformula}
\widehat{\Omega_{\lambda, \BB R}} (\phi) =
\int_{\g g_{\BB R}} F_{\lambda} \phi, \qquad
F_{\lambda}(X) = \sum_{x \in {\cal X}^X}
\frac {d_{X, x} e^{\langle X, \lambda_x \rangle}}
{\alpha_{x,1}(X) \dots \alpha_{x,n}(X)},
\end{equation}
where $d_{X, x}$'s are certain integers (to be described below)
which are local invariants of the sheaf ${\cal F}_{\lambda}$
given by formula (5.25b) in \cite{SchV2} and
$\alpha_{x,1}, \dots, \alpha_{x,n} \in \g t(X)^*$
are the roots of $\g t(X)$ such that, letting
$\g g^{\alpha_{x,1}}, \dots, \g g^{\alpha_{x,n}} \subset \g g$
be the corresponding root spaces and $\g b_x$ be the Borel subalgebra
of $\g g$ corresponding to $x$,
$$
\g b_x =
\g t(X) \oplus \g g^{\alpha_{x,1}} \oplus \dots \oplus \g g^{\alpha_{x,n}}
$$
as vector spaces.
Thus we obtain Harish-Chandra's formula for the invariant eigendistribution
$\widehat{\Omega_{\lambda, \BB R}}$.

It remains to specify the integer multiplicities $d_{X,x}$'s.
Choose a positive root system $\Phi_x^+$ on $\g t(X)$ so that, letting
$$
\g n_x^+ = \bigoplus_{\alpha \in \Phi_x^+} \g g^{\alpha},
$$
we have:
$$
\g g = \g b_x \oplus \g n_x^+
$$
as vector spaces.
Then we can select two subsets $\Phi_x', \Phi_x'' \subset \Phi_x^+$
with the following properties:
\begin{eqnarray*}
\text{a)} \quad &\text{for $\alpha \in \Phi_x^+$,
if $\operatorname{Re} (\alpha(X)) \ne 0$, then}  \\
&\alpha \in \Phi_x' \: \Longleftrightarrow \: \operatorname{Re} (\alpha(X)) < 0
\quad \text{and} \quad
\alpha \in \Phi_x'' \: \Longleftrightarrow \:
\operatorname{Re}(\alpha(X)) > 0;  \\
\text{b)} \quad &\text{for $\alpha_1, \alpha_2  \in \Phi_x^+$,
if $\alpha_1 + \alpha_2 \in \Phi_x^+$, then}  \\
&\alpha_1, \alpha_2 \in \Phi_x' \: \Longrightarrow \:
\alpha_1 + \alpha_2 \in \Phi_x', \quad
\alpha_1, \alpha_2 \in \Phi_x'' \: \Longrightarrow \:
\alpha_1 + \alpha_2 \in \Phi_x''.
\end{eqnarray*}
Such subsets  $\Phi_x', \Phi_x''$ always exist: for example,
the subsets defined by a) without the restriction
$\operatorname{Re}(\alpha(X)) \ne 0$ satisfy also b).
Because of b),
$$
\g n_x' = \bigoplus_{\alpha \in \Phi_x'} \g g^{\alpha}, \qquad
\g n_x'' = \bigoplus_{\alpha \in \Phi_x''} \g g^{\alpha}
$$
are subalgebras of $\g n_x^+$.
We also define
$$
N_x^+ = \exp (\g n_x^+) \cdot x, \qquad
N_x' = \exp (\g n_x') \cdot x, \qquad
N_x'' = \exp (\g n_x'') \cdot x.
$$
Then $N_x^+ \subset {\cal X}$ is an open Schubert cell which contains
$N_x'$ and $N_x''$ as affine linear subspaces.
Finally, the coefficient
\begin{equation}  \label{d_X,x}
d_{X,x} =
\chi \bigl( {\cal H}^{\bullet}_{N_x''} (\BB D {\cal F}_{\lambda})_x \bigr)
= \chi
\bigl({\cal H}^{\bullet}_{\{x\}} (\BB D {\cal F}_{\lambda} |_{N_x''}) \bigr)
\end{equation}
is the Euler characteristic, where $\BB D {\cal F}_{\lambda}$ denotes
the Verdier dual sheaf to ${\cal F}_{\lambda}$.
This multiplicity $d_{X,x}$ is exactly the local contribution of $x$ to
the Lefschetz fixed point formula generalized to sheaf cohomology by
M.~Goresky and R.~MacPherson \cite{GM}.

\begin{ex}  {\em
Suppose that the group $G_{\BB R}$ is compact. Then we can choose
$U_{\BB R} = G_{\BB R}$, in which case the twisted moment map
$\mu_{\lambda}: T^*{\cal X} \, \tilde \to \, \Omega_{\lambda}$
restricts to a $G_{\BB R}$-equivariant diffeomorphism
$\mu_{\lambda}: {\cal X} \, \tilde \to \, |\Omega_{\lambda,\BB R}|$.
We can choose a positive system of roots of $\g t$ so that
this restriction is orientation-preserving.
Then $C_{\lambda} = {\cal X}$, we can take
${\cal F}_{\lambda}= \BB C_{\cal X}$ -- the constant sheaf on ${\cal X}$,
all the coefficients $d_{X, x} =1$, and the formula (\ref{charformula})
coincides with the result of application of the classical fixed point
integral localization formula for equivariant forms
as described in Section 7.5 in \cite{BGV}.
}\end{ex}

%\separate

\begin{section}
{General Properties of the Cycle $C_{\lambda}$}
\end{section}  

In this section we describe some general properties of $C_{\lambda}$.
We begin with the following lemma.

\separate

The definition of the twisted moment map $\mu_{\lambda}$ involves choices
of a real compact form $U_{\BB R}$ and an identification between the
Cartan subalgebra $\g t \subset \g g$ containing $I^{-1}(i\lambda)$
and the universal Cartan algebra $\g h$.
We fix an isomorphism $\g t \simeq \g h$ and
consider two different compact real forms $U_{\BB R}, U_{\BB R}' \subset G$.
We denote by $\mu_{\lambda, U_{\BB R}}$
(respectively $\mu_{\lambda, U_{\BB R}'}$)
and $C_{\lambda, U_{\BB R}}$ (respectively $C_{\lambda, U_{\BB R}'}$)
the twisted moment map and the corresponding cycle
defined relatively to $U_{\BB R}$ (respectively $U_{\BB R}'$).

\begin{lem}
Suppose that $U_{\BB R}' = g U_{\BB R} g^{-1}$ for some $g \in G_{\BB R}$.
Then the cycles $C_{\lambda, U_{\BB R}}$ and $C_{\lambda, U_{\BB R}'}$
are the same.
\end{lem}

\pf
One can see from the definition of the twisted moment map that the relation
between the two maps $\mu_{\lambda, U_{\BB R}}$ and
$\mu_{\lambda, U_{\BB R}'}$ is
$$
\mu_{\lambda, U_{\BB R}'} (g \cdot \zeta) =
g \cdot \mu_{\lambda, U_{\BB R}} (\zeta),
\qquad \zeta \in T^*{\cal X}.
$$
Therefore,
\begin{multline*}
C_{\lambda, U_{\BB R}'} = \lim_{s \to 0^+}
s \cdot \mu_{\lambda, U_{\BB R}'}^{-1}(\Omega_{\lambda, \BB R})
= \lim_{s \to 0^+}
s \cdot g \cdot \mu_{\lambda, U_{\BB R}}^{-1}
( g^{-1} \cdot \Omega_{\lambda, \BB R})  \\
= g \cdot \lim_{s \to 0^+}
s \cdot \mu_{\lambda, U_{\BB R}}^{-1}(\Omega_{\lambda, \BB R})
= g \cdot C_{\lambda, U_{\BB R}}
= C_{\lambda, U_{\BB R}}
\end{multline*}
by $G_{\BB R}$-invariance of the cycles
$\Omega_{\lambda, \BB R}$ in $\Omega_{\lambda}$ and
$C_{\lambda, U_{\BB R}}$ in $T^*{\cal X}$.
\qed

\begin{lem}  \label{equivariant}
The composition  map
$\pi \circ \mu_{\lambda}^{-1} : \Omega_{\lambda} \to {\cal X}$
is $G$-equivariant.
\end{lem}

\pf
Recall that ${\cal X}$ is the set of all Borel subalgebras
$\g b \subset \g g$, $\g h$ is the universal Cartan algebra of $\g g$,
and $\lambda \in \g h^*$.
Then the map $\pi \circ \mu_{\lambda}^{-1}$ sends $\lambda \in \g h^*$
into the only Borel algebra $\g b$ which contains $\g h$ and all the
negative root spaces.
If $\nu \in \Omega_{\lambda}$ is an arbitrary element,
then $\nu = g \cdot \lambda$ for some $g \in G$, and $g \g h g^{-1}$ is
naturally a universal Cartan algebra of $\g g$.
On the other hand, $g \g b g^{-1} = g \cdot \g b$ is the Borel algebra
containing $g \g h g^{-1}$ and its negative root spaces.
Therefore, $\pi \circ \mu_{\lambda}^{-1}(\nu) = g \cdot \g b$.
This proves that the composition  map
$\pi \circ \mu_{\lambda}^{-1} : \Omega_{\lambda} \to {\cal X}$
is $G$-equivariant.
\qed

\begin{cor}  \label{orbit}
The composition  map
$\pi \circ \mu_{\lambda}^{-1} : \Omega_{\lambda} \to {\cal X}$
sends $|\Omega_{\lambda, \BB R}|$ into a single $G_{\BB R}$-orbit
in ${\cal X}$.
\end{cor}

Let ${\cal O}_{\lambda}$ denote this $G_{\BB R}$-orbit in ${\cal X}$:
$$
{\cal O}_{\lambda} = \pi \circ \mu_{\lambda}^{-1} (|\Omega_{\lambda, \BB R}|).
$$

We pick an element $\nu \in |\Omega_{\lambda, \BB R}|$, and let
$x = \pi \circ \mu_{\lambda}^{-1} (\nu) \in {\cal X}$.
Recall that $\g b_x$ denotes the Borel algebra corresponding to $x$,
i.e. the Lie algebra of the group $B(x)$ defined by (\ref{B(x)}).
Let $\g n_x = [\g b_x, \g b_x]$ be the nilpotent subalgebra of $\g b_x$,
and let $\g n_{x,\BB R} = \g n_x \cap \g g_{\BB R}$.

\begin{lem}  \label{subspace}
The coadjoint orbit $|\Omega_{\lambda,\BB R}|$ contains the set
$\nu + iI(\g n_{x,\BB R}) \subset i \g g_{\BB R}^*$.
\end{lem}

\pf
Let $N_{x,\BB R} = \exp(\g n_{x,\BB R})$ be the subgroup of $G_{\BB R}$
generated by $\g n_{x,\BB R}$. It is enough to show that
$$
N_{x,\BB R} \cdot \nu \quad \supset \quad \nu + iI(\g n_{x,\BB R})
$$
as subsets of $i \g g_{\BB R}^*$.
Moreover, it is enough to show that
$$
N_{x,\BB R} \cdot I^{-1}(\nu)
\quad \supset \quad
I^{-1}(\nu) + i \g n_{x,\BB R}
$$
as subsets of $i \g g_{\BB R}$.

First we notice that if $n \in \g n_{x,\BB R}$, then
$$
\exp(n) \cdot I^{-1}(\nu) = Ad(\exp(n)) (I^{-1}(\nu)) = e^{ad(n)} (I^{-1}(\nu))
\quad \subset \quad I^{-1}(\nu) + i \g n_{x,\BB R}
$$
because $ad(n) (I^{-1}(\nu)) \in i \g n_{x,\BB R}$.

We recall that the Lie algebra $\g n_{x,\BB R}$
is nilpotent and introduce its ideals
\begin{multline*}
\g n_{x,\BB R}^0 = \g n_{x,\BB R}, \quad
\g n_{x,\BB R}^1 = [\g n_{x,\BB R}, \g n_{x,\BB R}^0], \quad
\g n_{x,\BB R}^2 = [\g n_{x,\BB R}, \g n_{x,\BB R}^1], \quad \dots,  \\
\g n_{x,\BB R}^k = [\g n_{x,\BB R}, \g n_{x,\BB R}^{k-1}], \quad
\g n_{x,\BB R}^{k+1} = [\g n_{x,\BB R}, \g n_{x,\BB R}^k] = 0.
\end{multline*}
Because the element $I^{-1}(\nu) \in i \g g_{\BB R}$ is regular semisimple
and the centralizer of a regular semisimple element is a Cartan algebra
and hence consists of semisimple elements, each of the maps
$$
\g n_{x,\BB R}^j \to i \g n_{x,\BB R}^j, \qquad
n \mapsto [n,I^{-1}(\nu)],
$$
$j=0,\dots,k$ is a vector space isomorphism.
In particular, each of the maps
$$
\g n_{x,\BB R}^j / \g n_{x,\BB R}^{j+1} \to
i ( \g n_{x,\BB R}^j / \g n_{x,\BB R}^{j+1} ), \qquad
n + \g n_{x,\BB R}^{j+1} \mapsto [n,I^{-1}(\nu)] + i \g n_{x,\BB R}^{j+1},
$$
$j=0,\dots,k$, is a vector space isomorphism.
It follows that each of the maps
\begin{multline*}
\g n_{x,\BB R}^j \to
I^{-1}(\nu)+ i ( \g n_{x,\BB R}^j / \g n_{x,\BB R}^{j+1} ),  \\
n \mapsto
e^{ad(n)}(I^{-1}(\nu)) + i \g n_{x,\BB R}^{j+1} =
I^{-1}(\nu) + [n,I^{-1}(\nu)] + i \g n_{x,\BB R}^{j+1},
\end{multline*}
$j=0,\dots,k$, is surjective.
Therefore, we get the following inclusions:
$$
N_{x,\BB R} \cdot I^{-1}(\nu) + i \g n_{x,\BB R}^{j+1}
\quad \supset \quad
I^{-1}(\nu) + i \g n_{x,\BB R}^j,
$$
$j=0,\dots,k$, which together imply
$$
N_{x,\BB R} \cdot I^{-1}(\nu)
\quad \supset \quad
I^{-1}(\nu) + i \g n_{x,\BB R}.
$$
This finishes our proof of the lemma. \qed

\begin{prop}  \label{fiber}
The set $\mu_{\lambda}^{-1} (|\Omega_{\lambda, \BB R}|) \cap T_x^*{\cal X}$
is precisely $\mu_{\lambda}^{-1} (\nu) + T^*_{{\cal O}_{\lambda}, x} {\cal X}$,
where
$$
T^*_{{\cal O}_{\lambda}, x} {\cal X} = \{ \zeta \in T^*_x {\cal X} ;\:
\operatorname{Re} (\zeta)|_{T{{\cal O}_{\lambda}}} =0 \}
$$
is the ${\cal C}^{\infty}$ conormal space to ${\cal O}_{\lambda}$
at the point $x$.
\end{prop}

\pf
Pick any element $\zeta \in T_x^*{\cal X}$ such that
$\mu_{\lambda}^{-1} (\nu) + \zeta \in
\mu_{\lambda}^{-1} (|\Omega_{\lambda, \BB R}|)$.
Then
\begin{multline*}
i \g g_{\BB R}^* \supset |\Omega_{\lambda, \BB R}| \ni
\mu_{\lambda}( \mu_{\lambda}^{-1} (\nu) + \zeta ) =
\lambda_x + \mu( \mu_{\lambda}^{-1} (\nu) + \zeta )  \\
= (\lambda_x + \mu( \mu_{\lambda}^{-1} (\nu)) + \mu(\zeta) =
\nu + \mu(\zeta).
\end{multline*}
This shows that $\mu(\zeta) \in i \g g_{\BB R}^*$ or, equivalently,
$\zeta \in T^*_{{\cal O}_{\lambda}, x} {\cal X}$,
which proves one inclusion:
$$
\mu_{\lambda}^{-1} (|\Omega_{\lambda, \BB R}|) \cap T_x^*{\cal X}
\quad \subset \quad
\mu_{\lambda}^{-1} (\nu) + T^*_{{\cal O}_{\lambda}, x} {\cal X}.
$$

To prove the other inclusion we use the result of Lemma \ref{subspace} that
$$
|\Omega_{\lambda,\BB R}| \quad \supset \quad \nu + iI(\g n_{x,\BB R}).
$$
Notice that
$$
\mu_{\lambda}^{-1} \bigl( \nu + iI(\g n_{x,\BB R}) \bigr)
\quad \subset \quad T^*_x {\cal X}
$$
and also that
$$
\mu_{\lambda}^{-1} \bigl( \nu + iI(\g n_{x,\BB R}) \bigr)
- \mu_{\lambda}^{-1} (\nu)
\quad = \quad
\mu^{-1} \bigl( iI(\g n_{x,\BB R}) \bigr) \cap T^*_x {\cal X}
$$
is a real vector subspace of $T^*_x {\cal X}$.
The dimension of this subspace equals
\begin{multline*}
\dim_{\BB R} (\g n_{x,\BB R}) =
\dim_{\BB R}(\g g_{\BB R}) - \frac 12 \dim_{\BB C} (\g h) -
\dim_{\BB R} ({\cal O}_{\lambda})  \\
= \dim_{\BB R}({\cal X}) - \dim_{\BB R} ({\cal O}_{\lambda}) =
\dim_{\BB R} (T^*_{{\cal O}_{\lambda}, x} {\cal X}).
\end{multline*}
This forces
$$
\mu_{\lambda}^{-1} \bigl( \nu + iI(\g n_{x,\BB R}) \bigr)
- \mu_{\lambda}^{-1} (\nu)
\quad = \quad
T^*_{{\cal O}_{\lambda}, x} {\cal X}
$$
which implies the proposition.
\qed

\begin{thm}  \label{main}
\quad
\begin{enumerate}
\item
The cycle $C_{\lambda}$ can be written as
$$
C_{\lambda} = T^*_{{\cal O}_{\lambda}} {\cal X} + c_{\lambda},
$$
where the conormal space $T^*_{{\cal O}_{\lambda}} {\cal X}$
is equipped with an appropriate orientation and regarded as a chain
in $T^*{\cal X}$ and $c_{\lambda}$ is a chain supported inside
$\pi^{-1} (\overline{{\cal O}_{\lambda}} \setminus {\cal O}_{\lambda})$.
In particular,
%For $x \in {\cal O}_{\lambda}$, the set $|C_{\lambda}| \cap T_x^*{\cal X}$
%is precisely $T^*_{{\cal O}_{\lambda}, x} {\cal X}$ --
%the ${\cal C}^{\infty}$ conormal space to ${\cal O}_{\lambda}$
%at the point $x$. Hence
$$
|C_{\lambda}| \cap \pi^{-1}({\cal O}_{\lambda})
\quad = \quad T^*_{{\cal O}_{\lambda}} {\cal X}.
$$
\item
The coefficients in the fixed point formula (\ref{charformula})
$d_{X, x}$'s given by the equation (\ref{d_X,x})
are zero unless $x$ lies in the closure $\overline{{\cal O}_{\lambda}}$.
If $x \in {\cal O}_{\lambda}$, then $d_{X, x} = \pm 1$.
\end{enumerate}
\end{thm}

\pf
Part {\it (1)} follows form the definition of $C_{\lambda}$
(equation (\ref{C})), Corollary \ref{orbit} and Proposition \ref{fiber}.

Part {\it (1)} implies the first assertion of part {\it (2)}.
Let $\BB C_{{\cal O}_{\lambda}}$ denote the constant sheaf on
${\cal O}_{\lambda}$, and let
$j_{{\cal O}_{\lambda} \hookrightarrow {\cal X}}:
{\cal O}_{\lambda} \hookrightarrow {\cal X}$ denote the inclusion map.
Since the coefficients $d_{X,x}$'s described by the equation (\ref{d_X,x})
are determined by {\em local} properties of the cycle $C_{\lambda}$ at $x$,
as far as the second assertion of part {\it (2)} is concerned,
we can assume that
$$
C_{\lambda} = \pm Ch \bigl(
(R j_{{\cal O}_{\lambda} \hookrightarrow {\cal X}})_*
\BB C_{{\cal O}_{\lambda}} \bigr).
$$
(The choice of sign depends on the orientation of 
$T^*_{{\cal O}_{\lambda}} {\cal X}$.)
Then
$$
Ch(\BB D {\cal F}_{\lambda}) = \pm Ch \bigl(
(j_{{\cal O}_{\lambda} \hookrightarrow {\cal X}})_!
\BB C_{{\cal O}_{\lambda}} \bigr),
$$
and, for the purpose of computing the multiplicities $d_{X,x}$'s,
we may assume that $\BB D {\cal F}_{\lambda}$ is either
$$
(j_{{\cal O}_{\lambda} \hookrightarrow {\cal X}})_!
\BB C_{{\cal O}_{\lambda}}
\qquad \text{or} \qquad
(j_{{\cal O}_{\lambda} \hookrightarrow {\cal X}})_!
\BB C_{{\cal O}_{\lambda}}[1],
$$
where $[1]$ denotes a shift in degree by 1, so that
$$
Ch \bigl( (j_{{\cal O}_{\lambda} \hookrightarrow {\cal X}})_!
\BB C_{{\cal O}_{\lambda}} \bigr) =
- Ch \bigl( (j_{{\cal O}_{\lambda} \hookrightarrow {\cal X}})_!
\BB C_{{\cal O}_{\lambda}}[1] \bigr).
$$
Then the result follows immediately from (\ref{d_X,x}).
\qed

\begin{rem}
Part {\it (1)} of Theorem \ref{main} is precisely the assertion made
by W.~Rossmann in his proof of Theorem 1.4(b) in \cite{R3}.
(The proof of this statement was omitted because it was not the point
of that paper.)
\end{rem}

%\separate

\begin{section}
{The Maximally Split Cartan Case}  \label{splitsection}
\end{section}  

In this section we will identify the cycle $C_{\lambda}$
as a characteristic cycle of some sheaf ${\cal F}_{\lambda}$
when $I^{-1} (i\lambda)$ lies in a maximally split Cartan subalgebra
of $\g g_{\BB R}$.
Hence we will compute the function $F_{\lambda}$ in the formula
(\ref{F}) for the Fourier transform of $\Omega_{\lambda,\BB R}$.

\separate

We will use the {\em Iwasawa decomposition} of $G_{\BB R}$
(see \cite{Kn} for reference).
Let $U_{\BB R} \subset G$ be a compact real form, and let
$\g u_{\BB R} \subset \g g$ be its Lie algebra.
Let $\sigma$ and $\tau: \g g \to \g g$ denote the complex
conjugations with respect to $\g g_{\BB R}$ and $\g u_{\BB R}$
respectively. There is a choice of $U_{\BB R}$ such that
the involutions $\tau$ and $\sigma$ commute.
From now on we fix this $U_{\BB R}$ and let
$$
\theta = \sigma \tau = \tau \sigma: \g g \to \g g
$$
be the {\em Cartan involution} ($\theta^2=\operatorname{Id}_{\g g}$).
This involution $\theta$ preserves $\g g_{\BB R}$.
We define $\g k$ and $\g p$ as the $(+1)$ eigenspace
and the $(-1)$ eigenspace of $\theta: \g g \to \g g$
respectively, and let
$$
\g k_{\BB R} = \g k \cap \g g_{\BB R}, \qquad
\g p_{\BB R} = \g p \cap \g g_{\BB R}.
$$
Then
$$
\g g = \g k \oplus \g p, \qquad
\g g_{\BB R} = \g k_{\BB R} \oplus \g p_{\BB R}, \qquad
\g u_{\BB R} = \g k_{\BB R} \oplus i \g p_{\BB R}.
$$

We define $K_{\BB R}$ to be the identity component of
$U_{\BB R} \cap G_{\BB R}$.
Then $K_{\BB R}$ is a maximal compact subgroup of $G_{\BB R}$
with Lie algebra $\g k_{\BB R}$.
Let $\g a_{\BB R} \subset \g p_{\BB R}$ be a maximal abelian subspace,
let $\g m_{\BB R} = {\cal Z}_{\g k_{\BB R}}(\g a_{\BB R})$ be the
centralizer of $\g a_{\BB R}$ in $\g k_{\BB R}$, and let
$\g a = \g a_{\BB R} \otimes \BB C \subset \g g$ and
$\g m = \g m_{\BB R} \otimes \BB C \subset \g g$ be the
complexifications of $\g a_{\BB R}$ and $\g m_{\BB R}$.
Then we obtain a root space decomposition of $\g g$ with respect
to $\g a$:
$$
\g g = \g m \oplus \g a \oplus
\bigl( \bigoplus_{\alpha \in \Phi} \g g^{\alpha} \bigr),
$$
where $\Phi \subset \g a^* \setminus \{0\}$ is the reduction
of relative roots of the pair $(\g g, \g a)$.
We select a positive root system $\Phi^+ \subset \Phi$, and let
$$
\g n_{\BB R} = \g g_{\BB R} \cap \bigoplus_{\alpha \in \Phi^+} \g g^{-\alpha},
\qquad N_{\BB R} = \exp (\g n_{\BB R}) \subset G_{\BB R},
\qquad A_{\BB R} = \exp (\g a_{\BB R}) \subset G_{\BB R}.
$$
Then $G_{\BB R} = K_{\BB R} A_{\BB R} N_{\BB R}$
is the Iwasawa decomposition of $G_{\BB R}$.

\separate

We let $\g t_{\BB R} \subset \g m_{\BB R}$ be a maximal abelian
subalgebra consisting of semisimple elements,
then $\g h_{\BB R} = \g t_{\BB R} \oplus \g a_{\BB R}$
is a maximally split Cartan subalgebra of $\g g_{\BB R}$,
and its complexification $\g h = \g h_{\BB R} \otimes \BB C$
is a Cartan subalgebra of $\g g$.
We equip $\g h$ with a positive root system so that $\g n_{\BB R}$
lies in the vector space direct sum of the negative root spaces.
This way $\g h$ becomes a universal Cartan subalgebra.

In this section we compute the Fourier transform
$\widehat{\Omega_{\lambda, \BB R}}$ when
$\lambda \in i(\g g_{\BB R}^*)' \cap i I(\g h_{\BB R})$
is a regular semisimple element dual to the maximally split Cartan
subalgebra $\g h_{\BB R}$.

Notice that by construction the groups $A_{\BB R}$ and $N_{\BB R}$
fix the element $\pi \circ \mu_{\lambda}^{-1} (\lambda) \in {\cal X}$
and so the $G_{\BB R}$-orbit
$$
{\cal O}_{\lambda}
= G_{\BB R} \cdot \bigl( \pi \circ \mu_{\lambda}^{-1} (\lambda) \bigr)
= K_{\BB R} A_{\BB R} N_{\BB R} \cdot \bigl(
\pi \circ \mu_{\lambda}^{-1} (\lambda) \bigr)
= K_{\BB R} \cdot \bigl( \pi \circ \mu_{\lambda}^{-1} (\lambda) \bigr)
$$
is closed.
Then the first part of Theorem \ref{main} implies that the cycle
$C_{\lambda}$ is just the conormal space $T^*_{{\cal O}_{\lambda}}{\cal X}$
equipped with some orientation.

As in the proof of Theorem \ref{main},
let $\BB C_{{\cal O}_{\lambda}}$ denote the constant sheaf on
${\cal O}_{\lambda}$, and let
$j_{{\cal O}_{\lambda} \hookrightarrow {\cal X}}:
{\cal O}_{\lambda} \hookrightarrow {\cal X}$ denote the inclusion map.
Then
$$
C_{\lambda} = \pm Ch ({\cal F}_{\lambda}),
\qquad
{\cal F}_{\lambda} =
(j_{{\cal O}_{\lambda} \hookrightarrow {\cal X}})_*
\BB C_{{\cal O}_{\lambda}}.
$$

Substituting this sheaf into Harish-Chandra's formula (\ref{charformula})
and applying the second part of Theorem \ref{main} we obtain the following
description of the function $F_{\lambda}$ in the equation (\ref{F}).
It is an $Ad(G_{\BB R})$-invariant, locally $L^1$ function
on $\g g_{\BB R}$ which is represented by a real analytic function on
$\g g_{\BB R}'$.
If $X \in \g g_{\BB R}'$, then $F_{\lambda}(X)=0$ unless
$X$ is $G_{\BB R}$-conjugate to an element of $\g h_{\BB R}$.
And if $X \in \g h_{\BB R} \cap \g g_{\BB R}'$, then
$$
F_{\lambda}(X) = \sum_{x \in {\cal X}^X \cap {\cal O}_{\lambda}}
\frac {\pm e^{\langle X, \lambda_x \rangle}}
{\alpha_{x,1} (X) \dots \alpha_{x,n} (X)},
$$
where
${\cal X}^X$ denotes the set of zeroes of the vector field on ${\cal X}$
generated by $X$,
$\alpha_{x,1}, \dots, \alpha_{x,n} \in \g h^*$ are the roots of $\g h$
such that, letting
$\g g^{\alpha_{x,1}}, \dots, \g g^{\alpha_{x,n}} \subset \g g$
be the corresponding root spaces and $\g b_x$ be the Borel
subalgebra of $\g g$ corresponding to $x$,
$$
\g b_x = \g h \oplus
\g g^{\alpha_{x,1}} \oplus \dots \oplus \g g^{\alpha_{x,n}}
$$
as vector spaces.

%\separate

\end{document}